\documentclass{article}

\usepackage{amssymb}
\usepackage{amsmath}

\usepackage[dvips]{graphicx}

\begin{document}


\begin{center}
\textbf{On the continuous dual of the sequence space $ bv $ }
\end{center}
 
\noindent \textbf{}

\begin{center}
M. El Azhari 
\end{center}

\noindent \textbf{}

\noindent \textbf{Abstract.} Imaninezhad and Miri introduced the sequence space $ d_{\infty} $ in order to characterize the continuous dual of the sequence space $ bv. $ We show by a counterexample that this claimed characterization is false.
 
\noindent \textbf{}

\noindent \textbf{Mathematics Subject Classification 2010:} 46B10, 46B45.

\noindent \textbf{} 

\noindent \textbf{Keywords:} continuous dual, sequence space.

\noindent \textbf{}

\noindent \textbf{}

\noindent \textbf{1. Preliminaries }

\noindent \textbf{}                                                                        

\noindent \textbf{} Let $ \omega $ denote the vector space of all complex sequences, where addition and scalar multiplication are defined pointwise. For each $ k\geqslant 1, $ let $ e^{(k)}=(e_{n}^{(k)})_{n\geqslant 1} $ be the complex sequence defined by $e_{k}^{(k)}=1  $ and $ e_{n}^{(k)}=0 $ for $ n\neq k. $ A sequence space is a vector subspace of $ \omega $ including the set $ \{e^{(k)}:\:k\geqslant 1\}. $ By $ l_{\infty},\:c,$ and $ l_{q}\:(1\leqslant q<\infty), $ we denote the sequence spaces of all bounded, convergent and absolutely q-summable sequences, respectively. Also the inclusions 
$ l_{q} \subset c\subset l_{\infty}$ are strict. $ c $ and $ l_{\infty} $ are Banach spaces with the supremum norm $ \lVert x\rVert_{\infty}=\sup\{\lvert x_{k}\rvert: k\geqslant 1\}.\:\: l_{q} $ is a Banach space with the norm
 $ \lVert x\rVert_{q}=(\sum_{k=1}^{\infty}\lvert x_{k}\rvert^{q})^{\frac{1}{q}}. $
Let $ bv_{p}\:(1\leqslant p<\infty) $ denote the sequence space of all sequences of bounded variation defined, in [2], by $ bv_{p}=\{ x=(x_{k})_{k\geqslant 1}\in\omega:\:\sum_{k=1}^{\infty}\lvert x_{k}-x_{k-1}\rvert^{p}<\infty\} $ with $ x_{0}=0.\:\: bv_{p} $ is a Banach space with the norm $ \lVert x\rVert_{bv_{p}}=(\sum_{k=1}^{\infty}\lvert x_{k}-x_{k-1}\rvert^{p})^{\frac{1}{p}}. $ We will write $ bv $ instead of $ bv_{1}. $ If we define the sequence $ b^{(k)}=(b_{n}^{(k)})_{n\geqslant 1} $ of elements of the space $ bv_{p} $ for every $ k\geqslant 1 $
by $b_{n}^{(k)}=0  $ if $ n<k $ and $b_{n}^{(k)}=1  $ if $ n\geqslant k, $ then the sequence $ (b^{(k)})_{k\geqslant 1} $ is a Schauder basis [2] for $ bv_{p} $ and every $ x\in bv_{p} $ has a unique representation of the form 
$ x=\sum_{k=1}^{\infty} (x_{k}-x_{k-1})b^{(k)} $ with $ x_{0}=0. $
Let $ \Delta_{\omega} $ denote the difference operator on $ \omega $ defined by $ \Delta_{\omega}x=(x_{k}-x_{k-1})_{k\geqslant 1} $ with $ x_{0}=0 $ for every $ x=(x_{k})_{k\geqslant 1}\in\omega. $ For each $ 1\leqslant p<\infty, $ the map $ \Delta:bv_{p}\rightarrow l_{p},\:\Delta x=\Delta_{\omega}x, $ is an isometric linear isomorphism [2].
If $ (E,\lVert .\rVert) $ is a normed sequence space, then $ E^{\ast} $ denotes the continuous dual of $ E $ with the norm defined by $ \lVert f\rVert=\sup\{\lvert f(x)\rvert:\:\lVert x\rVert=1\} $ for all $ f\in E^{\ast}. $ 
It is well-known that $ l_{1}^{\ast}\cong l_{\infty} $ and $l_{p}^{\ast}\cong l_{q}\:(1<p,q<\infty)$ with $  \frac{1}{p}+\frac{1}{q}=1$                                                                                                                                                                                                                                                                                                                                                                                                                                                                                                                                                                                                                                                                                                                                                                                                                                                                                                                                                                                                                                                                                                                                                                                                                                                                                                                                                                                                                                                                                                                                                                                                                                                     

\noindent \textbf{}

\noindent \textbf{2. The sequence space $ d_{\infty}$} 

\noindent \textbf{}

\noindent \textbf{Proposition 2.1.} A linear functional $ f $ on $ bv $ is continuous if and only if there exists $ a=(a_{k})_{k\geqslant 1}\in l_{\infty} $ such that $ f(x)=\sum_{k=1}^{\infty}a_{k}(\Delta x)_{k} $ for all $ x\in bv. $ Furthermore, $ f(b^{(i)})=a_{i} $ and $ f(e^{(i)})=a_{i}-a_{i+1} $ for every $ i\geqslant 1. $

\noindent \textbf{}
 
\noindent \textbf {Proof.} The map $ \Delta:bv\rightarrow l_{1},\Delta x=(x_{k}-x_{k-1})_{k\geqslant 1} $ with $ x_{0}=0, $ is an isometric linear isomorphism [2]. Let $ f $ be a linear functional on $ bv. $ If $ f $ is continuous, then $ f\circ\Delta^{-1} $ is a continuous linear functional on $ l_{1}, $ so there exists $a=(a_{k})_{k\geqslant 1}\in l_{\infty} $ such that $ f\circ\Delta^{-1}(y)=\sum_{k=1}^{\infty}a_{k}y_{k} $ for all $ y\in l_{1}.$ Hence $ f(x)=(f\circ\Delta^{-1})(\Delta x)=\sum_{k=1}^{\infty}a_{k}(\Delta x)_{k} $ for all $ x\in bv. $ \\
Conversely, if $ f(x)=\sum_{k=1}^{\infty}a_{k}(\Delta x)_{k} $ for all $ x\in bv $ and some $ (a_{k})_{k\geqslant 1}\in l_{\infty}, $ then
$\vert f(x)\vert \leqslant\sum_{k=1}^{\infty}\vert a_{k}\vert\vert(\Delta x)_{k}\vert \leqslant \Vert a\Vert_{\infty}\sum_{k=1}^{\infty}\vert(\Delta x)_{k}\vert=\Vert a\Vert_{\infty}\Vert x\Vert_{bv}. $ Therefore $ f $ is a continuous linear functional on $ bv. $
Let $ i\geqslant 1, f(b^{(i)})=\sum_{k=1}^{\infty}a_{k}(\Delta b^{(i)})_{k}=a_{i} $ and $f(e^{(i)})=\sum_{k=1}^{\infty}a_{k}(\Delta e^{(i)})_{k}=a_{i}-a_{i+1}.$

\noindent \textbf{}

\noindent \textbf{Proposition 2.2.} $ bv^{\ast} $ is isometrically isomorphic to $ l_{\infty}. $
 
\noindent \textbf{}

\noindent \textbf{Proof.} Define $ \varphi:bv^{\ast}\rightarrow l_{\infty},\:\varphi(f)=(f(b^{(k)}))_{k\geqslant 1},\:\varphi  $ is a surjective linear map by Proposition 2.1, and $ \varphi $ is injective since $ (b^{(k)})_{k\geqslant 1} $ is a Schauder basis for $ bv. $ Let $ f\in bv^{\ast} $ and $ x\in bv,\: \lvert f(x)\rvert=\lvert f(\sum_{k=1}^{\infty}(\Delta x)_{k}b^{(k)})\rvert = \lvert \sum_{k=1}^{\infty}(\Delta x)_{k}f(b^{(k)})\rvert \leqslant \sum_{k=1}^{\infty}\lvert(\Delta x)_{k}\rvert \lvert f(b^{(k)})\rvert \leqslant\sup_{k\geqslant 1}\lvert f(b^{(k)})\rvert \sum_{k=1}^{\infty}\lvert(\Delta x)_{k}\rvert = \lVert\varphi(f)\rVert_{\infty} \lVert x\rVert_{bv}, $ then $ \lVert f\rVert \leqslant \lVert\varphi(f)\rVert_{\infty}. $ On the other hand, $ \lvert f(b^{(k)})\rvert \leqslant\lVert f\rVert\lVert b^{(k)}\rVert_{bv}=\lVert f\rVert $ since $ \lVert b^{(k)}\rVert_{bv}=1 $ for all $ k\geqslant 1, $ then $ \lVert\varphi(f)\rVert_{\infty} =\sup_{k\geqslant 1}\lvert f(b^{(k)})\rvert\leqslant\lVert f\rVert. $
 
\noindent \textbf{}
 
\noindent \textbf{} Imaninezhad and Miri [3] introduced the sequence space $ d_{\infty}=\{x=(x_{k})_{k\geqslant 1}\in\omega:\sup_{k\geqslant 1}\lvert\sum_{i=k}^{\infty}x_{i}\rvert <\infty\}, $ and claimed that $ T:bv^{\ast}\rightarrow d_{\infty},\: T(f)=(f(e^{(k)})_{k\geqslant 1}, $ is an isometric linear isomorphism [3, Theorem 3.3] (see also [4, Theorem 4.1]). Here we show that the set $ T(bv^{\ast}) $ is not included in $ d_{\infty}. $

\noindent \textbf{}
 
\noindent \textbf{Counterexample.} Let $ a=(a_{k})_{k\geqslant 1}\in l_{\infty}\smallsetminus c, $ and consider the map $ f(x)=\sum_{k=1}^{\infty}a_{k}(\Delta x)_{k} $ for all $ x\in bv. $ By Proposition 2.1, $ f\in bv^{\ast} $ and $ f(e^{(i)})= a_{i}-a_{i+1} $ for every $ i\geqslant 1. $ Let $ k\geqslant 1 $ and $ n\geqslant k,\:\sum_{i=k}^{n}f(e^{(i)})= \sum_{i=k}^{n}(a_{i}-a_{i+1})=a_{k}-a_{n+1}. $
If $ (f(e^{(k)}))_{k\geqslant 1}\in d_{\infty}, $ then $ \lim_{n\rightarrow\infty} \sum_{i=k}^{n}f(e^{(i)})=\sum_{i=k}^{\infty}f(e^{(i)})\in\mathbb{C}$ and consequently $ \lim_{n\rightarrow\infty}a_{n+1}=a_{k}-\sum_{i=k}^{\infty}f(e^{(i)}), $ which contradicts the fact that the sequence $ (a_{k})_{k\geqslant 1}  $ is not convergent.

\noindent \textbf{}

\noindent \textbf{Remark.} The above counterexample is also a counterexample to the equality $ \sum_{i=k}^{\infty}f(e^{(i)})=f(b^{(k)}) $ for all $ k\geqslant 1 $ and $ f\in bv^{\ast}, $ which is used in the proofs of [1, Theorem 2.3] and [3, Theorem 3.3].

\noindent \textbf{}

\noindent \textbf{3. The sequence space $ d_{q}\:(1<q<\infty) $}  
 
\noindent \textbf{}

\noindent \textbf{Proposition 3.1.} Let $ 1< p,q <\infty $ with $ \frac{1}{p}+\frac{1}{q}=1. $ A linear functional $ f $ on $ bv_{p} $ is continuous if and only if there exists $ a=(a_{k})_{k\geqslant 1}\in l_{q} $ such that $ f(x)=\sum_{k=1}^{\infty}a_{k}(\Delta x)_{k} $ for all $ x\in bv_{p}. $ Furthermore, $ f(b^{(i)})=a_{i} $ and $ f(e^{(i)})=a_{i}-a_{i+1} $ for every $ i\geqslant 1. $

\noindent \textbf{}

\noindent \textbf{Proof.} The map $ \Delta:bv_{p}\rightarrow l_{p},\Delta x=(x_{k}-x_{k-1})_{k\geqslant 1} $ with $ x_{0}=0, $ is an isometric linear isomorphism [2]. Let $ f $ be a linear functional on $ bv_{p}. $ If $ f $ is continuous, then $ f\circ\Delta^{-1} $ is a continuous linear functional on $ l_{p}, $ so there exists $a=(a_{k})_{k\geqslant 1}\in l_{q} $ such that $ f\circ\Delta^{-1}(y)=\sum_{k=1}^{\infty}a_{k}y_{k} $ for all $ y\in l_{p}.$ Hence $ f(x)=(f\circ\Delta^{-1})(\Delta x)=\sum_{k=1}^{\infty}a_{k}(\Delta x)_{k} $ for all $ x\in bv_{p}. $ \\
Conversely, if  $ f(x)=\sum_{k=1}^{\infty}a_{k}(\Delta x)_{k} $  for all $ x\in bv_{p} $  and some $ a=(a_{k})_{k\geqslant 1}\in l_{q}, $ then  $ \lvert f(x)\rvert\leqslant\sum_{k=1}^{\infty}\lvert a_{k}\rvert\lvert(\Delta x)_{k}\rvert \leqslant (\sum_{k=1}^{\infty}\lvert a_{k}\rvert^{q})^{\frac{1}{q}} (\sum_{k=1}^{\infty}\lvert (\Delta x)_{k}\rvert^{p})^{\frac{1}{p}}=\lVert a\rVert_{q} \lVert x\rVert_{bv_{p}} $ by H\"{o}lder inequality. Therefore $ f $ is a continuous linear functional on $ bv_{p}. $

\noindent \textbf{}

\noindent \textbf{Proposition 3.2.} Let $ 1< p,q <\infty $ with $ \frac{1}{p}+\frac{1}{q}=1. $ Then  $ bv_{p}^{\ast} $ is isometrically isomorphic to $ l_{q}. $   

\noindent \textbf{}

\noindent \textbf{Proof.} Define $ \psi:bv_{p}^{\ast}\rightarrow l_{q},\:\psi(f)=(f(b^{(k)}))_{k\geqslant 1},\:\psi  $ is a surjective linear map by Proposition 3.1, and $ \psi $ is injective since $ (b^{(k)})_{k\geqslant 1} $ is a Schauder basis for $ bv_{p}. $  Let $ f\in bv_{p}^{\ast} $ and $ x\in bv_{p},\: \lvert f(x)\rvert=\lvert f(\sum_{k=1}^{\infty}(\Delta x)_{k}b^{(k)})\rvert = \lvert \sum_{k=1}^{\infty}(\Delta x)_{k}f(b^{(k)})\rvert \leqslant \sum_{k=1}^{\infty}\lvert(\Delta x)_{k}\rvert \lvert f(b^{(k)})\rvert \leqslant 
(\sum_{k=1}^{\infty}\lvert(\Delta x)_{k}\rvert^{p})^{\frac{1}{p}}(\sum_{k=1}^{\infty}\lvert f(b^{(k)})\rvert^{q})^{\frac{1}{q}} = \lVert\psi(f)\rVert_{q} \lVert x\rVert_{bv_{p}} $ by H\"{o}lder inequality, then  $ \lVert f\rVert \leqslant \lVert\psi(f)\rVert_{q}. $
On the other hand, let $ f\in bv_{p}^{\ast} $ and define the sequence $ (x^{(n)})_{n\geqslant 1} $ of the space $ bv_{p} $ by $ (\Delta x^{(n)})_{k} = \lvert f(b^{(k)})\rvert^{q}f(b^{(k)})^{-1}$ if $ 1\leqslant k\leqslant n $ and $ f(b^{(k)})\neq 0,\:  (\Delta x^{(n)})_{k}=0 $
otherwise. Let $ n\geqslant 1,\: f(x^{(n)})= f(\sum_{k=1}^{\infty}(\Delta x^{(n)})_{k}b^{(k)}) = \sum_{k=1}^{\infty}(\Delta x^{(n)})_{k}f(b^{(k)})= \sum_{k=1}^{n}\lvert f(b^{(k)})\rvert^{q}. $\\
Then $ \sum_{k=1}^{n}\lvert f(b^{(k)})\rvert^{q}= f(x^{(n)})= \lvert f(x^{(n)})\rvert\leqslant\lVert f\rVert\lVert x^{(n)}\rVert_{bv_{p}}=\\
\lVert f\rVert (\sum_{k=1}^{n}\lvert f(b^{(k)})\rvert^{(q-1)p})^{\frac{1}{p}}= \lVert f\rVert (\sum_{k=1}^{n}\lvert f(b^{(k)})\rvert^{q})^{\frac{1}{p}} $ since $ (q-1)p=q. $
Therefore $ (\sum_{k=1}^{n}\lvert f(b^{(k)})\rvert^{q})^{\frac{1}{q}}= (\sum_{k=1}^{n}\lvert f(b^{(k)})\rvert^{q})^{1-\frac{1}{p}} \leqslant 
\lVert f\rVert $  for all $ n\geqslant 1. $ Letting $ n\rightarrow\infty, $ we get $ \lVert\psi (f)\rVert_{q}= (\sum_{k=1}^{\infty}\lvert f(b^{(k)})\rvert^{q})^{\frac{1}{q}} \leqslant\lVert f\rVert.$
 
\noindent \textbf{}

\noindent \textbf{} Akhmedov and Ba\c{s}ar [1] introduced the sequence space $ d_{q}=\{x=(x_{k})_{k\geqslant 1}\in\omega: \sum_{k=1}^{\infty} \lvert\sum_{i=k}^{\infty}x_{i}\rvert^{q} <\infty\}\:(1 < q <\infty), $ with the norm\\
$ \lVert x\rVert_{d_{q}}= (\sum_{k=1}^{\infty} \lvert\sum_{i=k}^{\infty}x_{i}\rvert^{q})^{\frac{1}{q}},$ and proved that the sequence space $ d_{q} $ is isometrically isomorphic to $ bv_{p}^{\ast}\: (1 < p <\infty) $ with $ \frac{1}{p}+\frac{1}{q}=1 $ (see [1, Theorem 2.3] or [3, Theorem 3.6]). Here we deduce this result as a consequence of Proposition 3.2.

\noindent \textbf{}

\noindent \textbf{Corollary 3.3.} Let $ 1< p,q <\infty $ with $ \frac{1}{p}+\frac{1}{q}=1. $ Then  $ bv_{p}^{\ast} $ is isometrically isomorphic to $ d_{q}. $ 

\noindent \textbf{}

\noindent \textbf{Proof.} By Proposition 3.2 and the fact that the map $ D:d_{q}\rightarrow l_{q},\: D((x_{k})_{k\geqslant 1})=(\sum_{i=k}^{\infty} x_{i})_{k\geqslant 1}, $ is an isometric linear isomorphism.
 
\noindent \textbf{}

\noindent \textbf{Remark.} By using the same above argument, we cannot show that $ bv^{\ast} $ is isometrically isomorphic to $ d_{\infty} $
since the linear map $ \delta:d_{\infty}\rightarrow l_{\infty},\:\delta((x_{k})_{k\geqslant 1})=(\sum_{i=k}^{\infty}x_{i})_{k\geqslant 1},\: $
which is injective and norm-preserving, is not surjective (we easily show that the element $ (1,1,1,...)\in l_{\infty} $ does not belong to the image of $ \delta $).

\noindent \textbf{}

\noindent \textbf{}

\noindent \textbf{}

\noindent \textbf{}

\begin{center}
REFERENCES
\end{center}
 
\noindent \textbf{} 
 
\noindent \textbf{} [1] Akhmedov A. M. and Ba\c{s}ar F., The fine spectra of the difference operator $ \Delta $ over the sequence space $ bv_{p}\:(1\leqslant p <\infty), $ Acta Math. Sin. Eng. Ser., 23(10)(2007), 1757-1768.

\noindent \textbf{} [2] Ba\c{s}ar F. and Altay B., On the space of sequences of p-bounded variation and related matrix mappings, Ukrainian Math. J., 55(1)(2003), 136-147.

\noindent \textbf{} [3] Imaninezhad M. and Miri M., The dual space of the sequence space $ bv_{p}\:(1\leqslant p <\infty), $ Acta Math. Univ. Comenianae, 79(1)(2010), 143-149.

\noindent \textbf{} [4] Imaninezhad M. and Miri M., The  continuous dual of the sequence space $ l_{p}(\Delta^{n}), $ Acta Math. Univ. Comenianae, 79(2)(2010), 273-280.

\noindent \textbf{} 

\noindent \textbf{} Ecole Normale Sup\'{e}rieure

\noindent \textbf{} Avenue Oued Akreuch

\noindent \textbf{} Takaddoum, BP 5118, Rabat

\noindent \textbf{} Morocco
 
\noindent \textbf{} 

\noindent \textbf{} E-mail:  mohammed.elazhari@yahoo.fr

\end{document}